\documentclass[12pt,leqno]{article}
\usepackage{amssymb,amsmath,amsfonts,amsbsy, xspace}

\newtheorem{theorem}{Theorem}
\newtheorem{lemma}{Lemma}

\newtheorem{definition}{Definition}
\newtheorem{corollary}{Corollary}

\newtheorem{remark}{Remark}
\newtheorem{example}{Example}

\newcommand{\asterisk}{\ensuremath{\textbf{*}\!}}

\newcommand{\N}{\ensuremath{\,\mathbb{N}}}
\newcommand{\Z}{\ensuremath{\,\mathbb{Z}}}

\newcommand{\Q}{\ensuremath{\,\mathbb{Q}}}
\newcommand{\R}{\ensuremath{\,\mathbb{R}}}

\newcommand{\K}{\ensuremath{\,\mathbb{K}}}

\newcommand{\Infinitesimal}[1]{\ensuremath{\,\mathcal{I}\!}{(\!#1)}}
\newcommand{\Finite}[1]{\ensuremath{\,\mathcal{F}\!}{(\!#1)}}
\newcommand{\InfLarge}[1]{\ensuremath{\,\mathcal{L}\!}{(\!#1)}}

\newcommand{\starR}{\ensuremath{\,\textbf{*}\mathbb{R}}}

\newcommand{\st}{\ensuremath{{\,\rm{st}}}}

\newcommand{\setValRing}[1]{\ensuremath{\mathcal{R}_v(\K) = \{ x \in
\K\mid\,v(x) \geq 0\,\}}}
\newcommand{\setValIdeal}[1]{\ensuremath{\mathcal{I}_v(\K) = \{ x \in
\K\mid\,v(x) > 0\,\}}}
\newcommand{\setValUnits}[1]{\ensuremath{\mathcal{U}_v(\K) = \{ x \in
\K\mid\,v(x) = 0\,\}}}

\newcommand{\Proof}{\noindent \textbf{Proof:\;}}

\begin{document}


\title{Back to Classics: Teaching Limits Through Infinitesimals\\
 }

\author{Todor D. Todorov (ttodorov@polymail.calpoly.edu)\\
                        Mathematics Department\\                
                        California Polytechnic State University\\
                        San Luis Obispo, California 93407, USA}
\date{}
\maketitle

\begin{abstract} The usual $\epsilon,\delta$-definition of the limit of a function (whether presented at a rigorous or an intuitive
level) requires a ``candidate $L$''  for the limit value. Thus, we have to start our first calculus course with ``guessing'' instead
of ``calculating''. In this paper we criticize the method of using calculators for the purpose of selecting candidates for $L$. We
suggest an alternative: a working formula for calculating the limit value L of a real function in terms of infinitesimals. Our
formula, if considered as a definition of limit, is equivalent to the usual $\epsilon,\delta$-definition but does not involve a
candidate $L$ for the limit value. As a result, the Calculus becomes  to ``calculate''  again as it was originally designed to do. 

\end{abstract}
\section{Introduction}\label{S: Introduction}

Let $f : X \to \R$ be a real function, where $X\subseteq\R$ and $r\in R$ be a non-trivial adherent (accumulation) point of $X$. Recall that a real number
$L$ is called the limit of  $f$  as  $x$ tends to $r$, in symbols,  $\lim_{x\to r}f(x) = L$, if 
\[
	(\forall\epsilon\in\R_+)(\exists\delta\in\R_+)\forall x \in X) ( 0 < | x - r | < \delta \implies | f(x) - L | < \epsilon).
\]
This is the so-called  $\epsilon,\delta$-definition of limit. It is sometimes attributed to Cauchy but it appears historically for
the first time in John Wallis's ``Arithmetica Infinitorum'' (The Arithmetic of Infinites) in 1655. In the next 250 years following
Wallis's work, this definition was rejected  and rediscovered many times until it was finally accepted by the mathematical community
in the beginning of the 20th century. We should mention that at that time most of the results in Calculus were already discovered
through infinitesimals. 
	Contemporary mathematicians might be puzzled by the fact that it took so long for the mathematical community to accept such a ``nice and rigorous''
definition, especially taking into account that its alternative in terms of infinitesimals was commonly viewed as ``certainly non-rigorous'' (although
practically efficient). We can detect at least three more obvious reasons for this amazing phenomena in the history of calculus:

		{\bf 1.} In the period ``from Leibniz to Weierstrass'', not only infinitesimals, but also real numbers, did not have a rigorous mathematical foundation.
Thus, although the $\epsilon,\delta$-definition makes perfect sense even in the framework of the rational numbers, this definition is completely
fruitless without the completeness of the reals. We can proudly declare now that the obstacles related to the completeness of the real numbers belong
to the past. Sadly, we can not be so proud about points 2 and 3 below.

		{\bf 2.} The  $\epsilon,\delta$-definition of limit is shockingly complicated due to the involvement of three non-commuting quantifiers $\forall, \exists,
\forall$. In the formulation for existence of a limit the quantifiers become four: $\forall, \exists, \forall, \forall$. As a result, a real analysis
course resembles a collection of exercises in mathematical logic rather than a rigorous version of calculus. The gap between the elementary calculus
and real analysis widens and some students understandably wonder whether these two branches of mathematics have anything in common. Considerable
efforts have been made to present the $\epsilon,\delta$-definition in a more digestible and human-like form mostly by using geometric language (L.
Gillman, R.H. McDowell~\cite{Gillman/Mc} and S. Lang~\cite{sLang}). It is almost a public secret, however, that most mathematicians think and do
research in terms of infinitesimals and use the $\epsilon,\delta$-definition of limits only to present the final version of their work in a socially
acceptable form. 

	While acknowledging the importance of the above two factors for the unusually slow and late acceptance of the
$\epsilon,\delta$-definition in the history of calculus, we would like to focus our attention on another unpleasant feature of the
$\epsilon,\delta$-definition :

	{\bf 3.} The $\epsilon,\delta$-definition of limit does not give any clue as to how to calculate the limit value $L$. 
	At least this is true in the framework of the class of {\bf all} functions (all polynomials or all rational functions, etc.), that is, 
functions with arbitrarily large derivatives. And this is exactly the situation in a typical first calculus course. 
Thus, we have to guess a reasonable value for $L$ and then prove or disprove the truthfulness of our guess with the help of the  $\epsilon,\delta$-definition.
If the graph of the function is known (or, for the contemporary mathematician, if it is already on our computer screen), then the value of $L$ can be
reasonably guessed. Notice, however, that in this case we hardly need the concept of limit. In a numerical analysis course we probably would restrict our
discussion to the class of functions with bounded derivatives (and given bound) and try to localize $L$ within a given interval. But how to find $L$ in a
calculus course for freshmen who, presumably, do not know either what ``limit'' is, nor what ``derivative'' is, let alone the class of ``functions with
bounded derivatives'' ?
	
		We sometimes are tempted to evaluate  $f$  for finitely many points different from  $r$ and try to guess the value of $L$ assuming that there exists some
pattern in the behavior of the function (J. Stewart~\cite{jStew}, p. 50-61). Suppose, for example, that $f(r + 10^{-10}) = 1.99999999999$. Maybe the correct
limit of  $f$  (as $x$  approaches $r$) is $L = 2$ ? This guess is rooted in the following two myths:

		{\bf Myth 1:} $10^{-10}$  is a ``small number''. More generally, ``there are numbers in \R\, which are small and others which are large''. For example,
$10^{10}$ is ``certainly a large number''. Indeed, we never use ``millimeters'' to measure the distance between two cities, nor do we count our annual income
in ``cents''. By changing the units of measurement (to suit our convenience) we always try to stay away from numbers such as $10^{-10}$ or $10^{10}$. The
usage of different units for measurement explains the origin of this myth which, of course, has nothing to do with mathematics. 
	
	{\bf Myth 2 :} This myth originates in our experience as high-school students. It says that : ``The integers are more likely to be the correct answers than
are the fractions''. Thus, we somehow prefer to believe that $L = 2$, not $L = 1.999999999998$, is the correct answer. 

		Unfortunately, the values of a function  $f$ at finitely many points do not determine uniquely the limit of the function. We have to use even stronger
language : The evaluation of a given function at finitely many points (different from the limit point $r$) is completely and totally irrelevant to both the
value of the limit $L$ and to the concept of limit in general. Here is the precise negative statement:  

\begin{lemma}\label{L: Do not Guess}{\em (Do not Guess !):} Let $f : X\to \R,\; X\subseteq\R$, be a real function and $r \in \R$ be a non-trivial adherent
(accumulation) point of $X$. Let $P(x_i, y_i), i = 1, 2, ..., m$,  be finitely many (distinct) points in the plane  $\R^2$ such that $x_i \in X, x_i
\not= r$ and also such that  $x_i = x_j$ implies $y_i = y_j$. Let $L \in \R$ be an arbitrarily chosen real number (or even  $\pm \infty$). Then there
exists a polynomial  $f : X \to \R$ (or, a rational function  $g$ ) such that $f(x_i) = y_i, i = 1, 2, ..., m$, and   $\lim_{x\to r}f(x) = L$ (or,   
$\lim_{x\to r}g(x) = \pm\infty$, respectively).  
\end{lemma}
\Proof : Choose a polynomial  $f$  of  degree $m$ 
and solve the linear system of equations $f(x_i) = y_i, i = 1, 2, ..., m, f(r) = L$, for the coefficients in  $f$. The system has always at least one solution
and we have  $\lim_{x\to r}f(x) = L$, since  $f$  is continuous at $r$. $\blacktriangle$

	In Section 2 we intend to show that if we use an arbitrarily chosen non-zero infinitesimal $dx$ (instead of the increment $10^{-10}$ in our earlier
example), then the value $f(r + dx)$ uniquely determines the limit $L$. In Section 3 we present several examples from calculus to demonstrate how our method
works in practice. For those readers who remain skeptical toward the practical merits of our approach, we remind them that all inventions of what we call
today ``calculus'' have historically been discovered by means of infinitesimals. The reader who is interested in the history of calculus  (C.H. Edwards,
Jr.~\cite{cEdwards} and A. Robinson~\cite{aRob66}, Chapter X) will certainly observe that our method for calculating limits - if applied to calculating
derivatives - is similar to the original Leibniz-Euler infinitesimal method. This explains the phrase: ``Back to Classics'' in the title of our article.
	
	This article is written for calculus teachers who are looking for alternatives to the conventional methods for teaching limits. We shall try to 
keep the exposition at the level of rigor in which complex numbers are defined: as ``expressions of the form $x + i y$'', where $x$ and $y$ are real numbers
and $i =\sqrt{-1}$. A more advanced (but still accessible) introduction to Infinitesimal Calculus is presented in Section~\ref{S: An Introduction to
Infinitesimal Calculus} of this paper, where the reader will find precise definitions and complete proofs. 

\section{Preliminaries: Standard Part Mapping and Hyperreal Numbers}

We introduce the concept of infinitesimal and study the basic properties of hyperreal (nonstandard) numbers. We also study the operation known as the
standard part mapping which is, in a sense, an algebraic counterpart of the concept of limit applied to numbers rather than to functions.  

\begin{definition}\label{D: Infinitesimals, Finite and Infinitely Large Numbers}{\em (Infinitesimals, Finite and Infinitely Large Numbers):}

	{\bf (i)} A number  $dx$   is called {\bf infinitesimal},  in symbols,  $dx \approx 0$,  if  
$| dx | < 1/n$  for all $ n \in\N$. If  $x - y$ is an infinitesimal, we say that  $x$  and $y$  are infinitely close, in symbols,  $x \approx y$. 

	{\bf (ii)} A number  $x$ is called {\bf finite}  if  $| x |\leq n$ for some $n \in \N$. 

	{\bf(iii)}  A number  $x$   is called {\bf infinitely large} if  $n < |x|$  for all $n\in\N$.
	
	{\bf(iv)} If  $x$  and $y \not= 0$ are finite numbers, then the numbers of the form $x/y$ are called hyperreal (or nonstandard) numbers. The set of all
hyperreal (nonstandard) numbers will be denoted by {\em \starR}. \par
	We denote by {\em \Infinitesimal{\starR}, \Finite{\starR} and \InfLarge{\starR}} the sets of thy
infinitesimal, finite and infinitely large numbers in {\em \starR}, respectively.
\end{definition}

	It is clear that all real numbers are finite and zero is the only infinitesimal in \R. Also \R\, has no infinitely large numbers.
	The following rules follow directly  from the above definition:

\begin{theorem}\label{T: Properties} {\em (Properties):}

		{\bf (i)} finite   $\pm$   finite  $=$ finite,\; finite  $\times$ finite $ =$ finite.

		{\bf (ii)} infinitesimal  $\pm$  infinitesimal $ =$ infinitesimal.

		{\bf (iii)} infinitesimal  $\times$ infinitesimal $ =$ infinitesimal. 

		{\bf (iv)} infinitesimal $\times$ real $=$ infinitesimal and, more generally, infinitesimal  $\times$  finite $=$ infinitesimal.

		{\bf (v)} positive infinitely large  $+$  positive infinitely large $=$ positive infinitely large.
 
		{\bf (vi)} positive infinitely large $\times$  positive infinitely large 	$=$ positive infinitely large. 
	
	 {\bf (vii)} 1/non-zero infinitesimal $=$ infinitely large number.
\end{theorem}

\noindent{\bf Warning:} Numbers of the forms: $\frac{\textrm{infinitesimal}}{\textrm{infinitesimal}}$, $\frac{\textrm{finite}}{\textrm{finite}}$, $\frac{\textrm{infinitely large}}{\textrm{infinitely large}}$, ``positive (negative) infinitely large  -  positive (negative) infinitely large'' are always well defined provided that the denominators (if any) are non-zero. However, they can be of any type: infinitesimal, finite (real) or infinitely large. For example, let $dx$ be a non-zero infinitesimal, in symbols, $dx \neq 0,\, dx \approx 0$. Then $dx^2/dx$  is infinitesimal, $dx/dx$ is a real number (and that is $1$), $dx/dx^2$ is infinitely large, both $(2+dx)/(2+dx)^2$  and $(2 + dx)^2/(2 + dx)$ are finite,  $dx^{-1}/dx^{-2}$ is infinitesimal and   $dx^{-2}/dx^{-1}$ is infinitely large. Finally, $1/dx - 1/dx^2$  is infinitely large, $(1/dx + dx) - 1/dx = dx$ is infinitesimal and $(1/dx + 5) - 1/dx = 5$  is finite (actually, real). 

\begin{remark} {\em The level of rigor of Definition~\ref{D: Infinitesimals, Finite and Infinitely Large Numbers} is similar to the level of rigor of the
definition of complex numbers as ``expressions of  the form $x + i y$, where $x$ and $y$ are real numbers and  $i = \sqrt{-1}$\, ''. It is clear that our
definition "hangs on" the existence of a non-zero infinitesimal (just as the existence of complex numbers ``hangs on'' the existence of\,  $i = \sqrt{-1}$).
On the other hand, it is clear that if there exists one non-zero infinitesimal $dx$, then there are infinitely many infinitesimals: $2dx, 3dx, 4dx, dx^2,
dx^3$, etc. are also infinitesimals.
}\end{remark}

\noindent{\bf Axiom 1:} There exists a non-zero infinitesimal $dx$, in symbols, $dx \neq0,\, dx\approx 0$. 

\begin{remark}\label{R: TwoPrinciples}{\em The above definition together with Axiom 1 is a folk-like version of the statement: {\bf Let \starR\,  be a proper
totally ordered field extension of  \R}. We are simply trying to avoid fancy terminology. Recall that every totally ordered proper field extension of \R\,
is a non-Archimedian field, hence,  it contains non-zero infinitesimals and infinitely large numbers. Strictly speaking, not any proper totally ordered
field extension
\starR\, of \R\, is adequate for the needs of Calculus.  We also need that \starR\, is a ``non-standard extension''  of \R\, which means that \starR\,
satisfies two additional axioms. The first axiom  (Transfer Principle) says, roughly speaking, that every function $f$  in $\R^d$ ($d\in\N$) has an
extension 
$\asterisk f$ in $\starR^d$ such that the mapping\,  $*$\, preserves the equivalence between the equations and inequalities in \R\,
and \starR,\, respectively (where the right and left hand sides of these equations and inequalities are considered as real functions and their
${\bf *}$-extensions, respectively). For example, we have
$-1 \leq \sin{x} \leq 1 \iff  x = x $ in \R. Hence, by the Transfer Principle, it follows  $-1 \leq {\asterisk}\sin{x} \leq 1 \iff  x = x $ in 
\starR. In other words, the range of\, $\asterisk\sin{x}$ is the set $\{\, y \in \starR \mid -1 \leq y\leq 1\,\}$. We should mention that the
concept of ``nonstandard extension of a field'' is in sharp contrast to the concept of an ``algebraic extension of a field'' (where the equivalence between
some polynomial equations is, by design, violated in the extended field). An important consequence of the Transfer Principle is that \starR\, is a {\bf real
closed field}, meaning that every polynomial equation of odd degree with coefficients in \starR\, has a solution in \starR. The second axiom (Saturation
Principle) is a sort of completeness which implies, in particular, that every nested sequence of open intervals in
\starR\, has a non-empty intersection.  For a more detailed exposition of nonstandard analysis by means of these two Principles we refer to the
Keisler's Calculus textbook and its companion, written for calculus instructors (H. J. Keisler~\cite{hjKeisE}-\cite{hjKeisF}). Both axioms (especially the
Transfer Principle) are so natural that only an experienced mathematician will realize that they are, actually, needed in Calculus. The situation is similar
to the role of the Axiom of Choice in Real Analysis. It is well known that Real Analysis can not survive without the Axiom of Choice but it is completely
possible to teach Real Analysis without even mentioning this axiom.  }\end{remark}

If the reader still does not feel comfortable with the definition of hyperreal numbers, presented above, he/she should refer (now or later) to
Section~\ref{S: An Introduction to Infinitesimal Calculus} at the end of this paper. We resume our discussion on teaching calculus. 
	It is clear that if $r$  is a real number and $dx$ is an infinitesimal, then $ r + dx$ is a finite number. Due to the completeness of \R,\, the reverse is
also true:

\begin{theorem}\label{T: Asymptotoc Expansion}{\em (Asymptotic Expansion):}  Every finite number $x$ in {\em \starR}\, has an asymptotic expansion of the
form $x = r + dx$  for some real number  $r \in \R$\, and some infinitesimal $dx \approx 0$. 
\end{theorem}

\Proof See Theorem~\ref{T: Standard Part Mapping} in the last section of this paper. $\blacktriangle$

\begin{remark}{\em (Completeness of  \R):} {\em H. J. Keisler~\cite{hjKeisF}, p.17-18) proved that the statement of the above theorem is, actually,
equivalent to the order completeness of \R\, (if \R\, is treated merely as a ``totally ordered field''). It is worth noticing that the completeness of
the real numbers in the form presented above appeared (treated as an obvious fact) in the early Leibniz-Euler Infinitesimal Calculus - 150 years
before Cauchy, Bolzano, Weierstrass and Cantor formulated the completeness of the reals in the forms known from the contemporary real analysis
textbooks. What does all this mean ? Well, perhaps the theory of the real numbers at the time of Leibniz and Euler was not so non-rigorous after all;
it only takes so long until we finally figure out how the creators of Calculus preferred to express the completeness of \R. 
}\end{remark}

	In addition to the above properties of the finite numbers we have the following uniqueness result:

\begin{lemma}\label{L: Uniqueness}{\em (Uniqueness):}  Let  $r\in \R$\, and  $dx \not= 0$. Then $r + dx = 0$ implies both $r = 0$  and  $dx = 0$. 
\end{lemma}

\Proof  $r + dx = 0\quad if{f}\quad  r = - dx$. Hence, r = 0, since the zero  0  is the only infinitesimal in \R.  $\blacktriangle$

	The above property justifies the following definition.

\begin{definition}\label{D: Standard Part Mapping}{\em (Standard Part Mapping):}  We define {\em $\st : \starR\to \R\,\cup\{\pm \infty\}$} by: 
	
{\bf (a)} $\st(r + dx) = r$  for  $r\in\R,\,  dx \approx 0$;

	{\bf (b)} $\st (1/dx) =\pm\infty,\, dx \approx 0$, for $dx > 0$ and $dx < 0$, respectively. 
\end{definition}

	In the case of finite numbers, the above several results can be summarized in the following corollary:

\begin{corollary}\label{C: Asymptotic Expansion of Finite Numbers}{\em (Asymptotic Expansion of Finite Numbers):} Every finite number 
{\em $x\in\Finite{\!\starR}$}  can be presented uniquely in the form $x  =  st(x)  +  dx$, where  $dx  =  x - st(x)$ is infinitesimal. We shall sometimes
refer to the above formula as an asymptotic expansion  of\,  $x$. 
\end{corollary}

	The next result follows immediately from the definition of "infinitesimal".

\begin{theorem} \label{T: Properties of st} {\em (Properties of \st):} Let $x$ and $y$ be finite numbers.  Then we have: 

	{\bf (i)}  	$x \approx y$\quad  $if{f}$\quad  $\st(x) = \st(y)$. In particular, $\st(dx) = 0$ for all infinitesimals  $dx$.

	{\bf (ii)}  $\st(r) = r$ for all  $r \in \R$.

	{\bf (iii)}  Let $x$ and $y$ be not infinitely close. Then $x < y\; if{f}\; \st(x) <  \st(y)$. Consequently, for arbitrary $x$ and $y$, 
``$x < y$ implies $\st(x)\leq \st(y)$'' and  ``$x \leq y$  implies $\st(x)\leq \st(y)$''. 

	{\bf (iv)} $\st(x \pm y) = \st(x) \pm \st(y)$.

	{\bf (v)}  $\st(x\,y) = \st(x)\,\st(y)$.

	{\bf (vi)}   $\st(x/y) = \st(x)/\st(y)$ whenever $\st(y) \not= 0$.
 
	{\bf (vii)}  $\st (x^n) = (\st(x))^n$  for all  $n \in \N$.

	{\bf (viii)} $\st(\sqrt[n]{x}) = \sqrt[n]{\st(x)},  n \in\N$, where the condition $x > 0$ (implying $st(x) \geq 0$ )  is required in the case of even
$n$.
\end{theorem}

\begin{remark} \label{R: Field, Ring, Ideal}{\em (Field, Ring, Ideal):} {\em In the usual algebraic terminology the results of the above theorem can be
summarized as follows: The set of hyperreal numbers \starR\, is a totally ordered non-Archimedean real closed field, the set of finite numbers
\Finite{\starR}\, is a convex ring without zero-divisors (a totally ordered integral domain), the set of infinitesimals  {\Infinitesimal{\starR}}\,  is
a convex maximal ideal in  \Finite{\starR}\,  and the factor space \Finite{\starR}\,F/\Infinitesimal{\starR}\, is isomorphic to \R\,
under
\st. 
}\end{remark} 
\begin{remark}\label{R: Extended Real Line}{\em (Extended Real Line):} {\em If the standard part mapping $\st$ acts on infinitely large numbers, the
result is either $\infty$ or $-\infty$. In these cases we have to perform the usual legal  and illegal  operations in the extended real line
$\R\cup\{\pm\infty\}$:

	{\bf (a)} Legal Operations:  All Operations in \R\, are legal. In addition, the following are also legal ($\epsilon$ is a positive real number). 
	\begin{align}\notag
&\infty +\infty =\infty,\;  - \infty - \infty = - \infty,\;   \pm\epsilon + \infty = \infty,\;   \pm\epsilon -\infty = -\infty,\\ \notag
&\infty \times\infty =\infty,\,  (-\infty)\times (-\infty) = \infty,\; (-\infty) \times \infty = - \infty, \\\notag 
&\pm \epsilon \times\infty = \pm\infty,\;   \epsilon\times (\pm\infty) = \pm\infty,\; -\epsilon \times (\pm\infty) =\mp\infty,  \\\notag 
&1/\pm\infty=0,\; \ln(\infty) = \infty,\;   \epsilon^{-\infty} =  0,\;   e^\infty  = \infty. \notag 	
\end{align}
	{\bf (b)} Illegal Operations include: 
\[
   		\infty - \infty,\quad\frac{1}{0},\quad \frac{0}{0},\quad\frac{\infty}{\infty},\quad 0 \times (\pm \infty),\quad 0^0.
\]
}\end{remark}
		
	Here are several exercises with standard part: 
\begin{example}\label{Ex: Exercises with st}{\em (Exercises with \st):} {\em In what follows  $r$  denotes a real number, and $dx$  and $dy$  are non-zero
infinitesimals.

	{\bf 1.} $\st(dx) = 0$. Similarly, $\st(3dx - 4dx^2 ) = 0$ since $3dx - 4dx^2$ is an infinitesimal.  

	{\bf 2.} $\st(\frac{-3 + dx - dx^2}{2 + 2dx\,dy + dy^3}) = \frac{\st(-3 +dx - dx^2)}{\st(2 + 2dx\,dy + dy^3)} = -3/2$.

	{\bf 3.} {\bf Incorrect:} $\st\left(\frac{\sqrt{3+dx}-\sqrt{3}}{dx}\right) = \frac{\st(\sqrt{3+dx}-\sqrt{3})}{\st(dx)} = 
\frac{\st(\sqrt{3+dx})-\st(\sqrt{3})}{\st(dx)}  =  \frac{0}{0}$, which is an illegal symbol. We disregard this calculation and try something else:
{\bf Correct:} 
\begin{align*}
   \st\left(\frac{\sqrt{3+dx}-\sqrt{3}}{dx}\right) 
   &= \st\left(\frac{(\sqrt{3+dx}-\sqrt{3})(\sqrt{3+dx}+\sqrt{3})}{dx\,(\sqrt{3+dx}+\sqrt{3})}\right)\\
   &= \st\left(\frac{3+dx-3}{dx\,(\sqrt{3+dx}+\sqrt{3})}\right) \\
   &= \st\left(\frac{dx}{dx\,(\sqrt{3+dx}+\sqrt{3})}\right)\\
   &= \st\left(\frac{1}{\sqrt{3+dx}+\sqrt{3}}\right)= \frac{\st(1)}{\st(\sqrt{3+dx}+\sqrt{3})} \\
   &= \frac{1}{\sqrt{\st(3+dx)}+\st(\sqrt{3}))}=\frac{1}{\sqrt{3}+\sqrt{3}} = \frac{1}{2\,\sqrt{3}}
\end{align*}
 , which is the correct answer. 
 
	{\bf 4.} {\bf Incorrect:} $\st\left(\frac{-3+dx}{2dx+dx^2}\right) = (\frac{\st(-3+dx)}{\st(2dx+dx^2)} =-3/0$, which is an illegal symbol in $\R\,\cup\{\pm
\infty\}$. As before, we have to disregard this calculation and try something else: {\bf Correct:} $\st\left(\frac{-3+dx}{2dx+dx^2}\right) =
\st\left(\frac{1}{dx}\,\frac{-3+dx}{2+dx}\right) =
\st\left(\frac{1}{dx}\right)\st\left(\frac{-3+dx}{2+dx}\right) =\pm\infty\times(-3/2) = \pm\infty$, depending upon whether  $dx$  is positive or
negative, respectively.
}\end{example}

\begin{remark}\label{R: Guide}{\em (Guide):} {\em The following guide might help us to decide ``what to do next'' when we calculate the standard part 
$\st(x)$ of a hyperreal number $x$: 

	{\bf (a)} Let $x$  be a finite number initially presented in the form $x = r + dx$ (or it can be easily presented in this form).
Then in order to calculate st(x) we have simply to apply the definition of \st, i.e. to "drop the infinitesimal term $dx", \st(x) = \st(r + dx) = r$.

	{\bf (b)} Let  $x$  be a hyperreal number of unknown type (it might be a finite number but not presented in the form $x = r + dx$, or it might be an
infinitely large number). To calculate $\st(x)$, we have to apply any of the properties of st, presented in Theorem~\ref{T: Properties of st}, and/or any of
the legal operations in the extended real line $\R\,\cup \{\pm\infty\}$ (Remark~\ref{R: Extended Real Line}). If the result of this calculation is a real
number  or a legal symbol in $\R\,\cup \{\pm\infty\}$ (part (a) in Remark~\ref{R: Extended Real Line}), then this is the correct answer for $\st(x)$. If at
some stage of our calculations we obtain an illegal symbol in $\R\,\cup \{\pm\infty\}$ (part (b) in Remark~\ref{R: Extended Real Line}), we should stop,
disregard the work done so far and start from the beginning a trying a different algebraic strategy.  }\end{remark}
\begin{remark}{\em (Is the Algebra Familiar ?):} {\em The reader should not be surprised that the ``algebra in the above examples sounds
familiar'' since  both \R\, and \starR\, are real closed fields and, as we know, all real closed fields (Archimedean or not) obey the
same laws of algebra.    }\end{remark}

\section{Limits Using Infinitesimals: Our Working Formulas}

	In this section we derive several formulas for calculating the limit value $L$ of a real function in terms of infinitesimals. Our working formulas -
if considered as definitions - are equivalent to the $\epsilon,\delta$-definitions of the corresponding limits in real analysis but they do not involve a
candidate
$L$ for the limit value. Thus, we can start teaching a calculus course by ``calculating'' instead of ``guessing and proving''. At the end of this section we
summarize the main features of our method. \\

\noindent{\bf Warning:} The level of the following exposition is slightly higher than it is appropriate for teaching in class. To make it more accessible we
recommend to following:

	{\bf (a)} All details connected with the domain of the function should be skipped. They might be discussed later in the course when (and if) necessary.

{\bf (b)} The question of the existence of limits (which we discuss below) should be left aside or discussed when (and if) this question arises
naturally in class. 

	{\bf (c)} The different types of limits should be presented in different sessions (not all at once as we have done below).

	{\bf (d)} In the beginning the emphasis should be on those limits which have immediate geometric applications: the vertical and horizontal
asymptotes. These limits are more important for the purpose of ``sketching the curve'' than the limits of the type $\lim_{x\to r} f(x) = L$, where
both  $r$  and  $L$ are real numbers (not $\pm\infty$). The latter, although of fundamental importance for calculus, have a more subtle meaning, mostly
to support the theory of ``continuity'' and ``derivatives.''  

\begin{definition} \label{D: Limit}{\em (Limit):} Let $f : X \to \R,\; X \subseteq\R$, be a real function and $(a, r)\cup (r, b)\subseteq X$
for some $a, b, r \in R,  a < r < b$. Suppose that   
\begin{equation}\label{E: Exists}
	\st[f(r + dx)] = \st[f(r + dy)]\quad \text{in}\quad \R\cup\{\pm\infty\},
\end{equation}
for all non-zero infinitesimals dx and dy. Then $\st[f(r + dx)] \in \R\cup\{\pm\infty\}$ is called the limit of   $f$  as  $x$ 
approaches $r$, in symbols,
\begin{equation}\label{E: WorkingFormula}
\lim_{x\to r} f(x) = \st[f(r + dx)],  
\end{equation}
where $d$x in the latter formula is an arbitrarily chosen non-zero infinitesimal. We shall refer to (\ref{E: WorkingFormula}) as our ``working formula''
for calculating limits through infinitesimals.
\end{definition}

\begin{remark}\label{R:  Existence of Limit}{\em (Existence of Limit):} {\em Notice that the condition (\ref{E: Exists}) guarantees the existence of the
corresponding limit value  $\st(f(r + dx)) = L$  and the independence of our working formulae on the choice of $dx$. We do not need to check the
condition (\ref{E: Exists}) before applying (\ref{E: WorkingFormula}). Rather, we should start with (\ref{E: WorkingFormula}) and when the calculations
are done and the value $L$ is obtained, we should check whether the value $L$ depends on the choice of
$dx$. If the answer is ``no'', then $L$ is the desired limit value. If the value of $L$ depends on the choice of $dx$ (say, $L$  might depend on the sign
of $d$x), then the limit $\lim_{x\to r} f(x)$ does not exist. As we already mentioned, it is preferable to skip the discussion of the existence of
limit  and focus on the working formula (\ref{E: WorkingFormula}).  
}\end{remark}

\begin{definition}\label{D: One Side Limits}{\em (One Side Limits):}
 
	{\bf (i)} Let $f : X \to \R,\; X\subseteq \R$, be a real function and $(r, r + \epsilon) \subseteq X$ for some $r\in\R$ and some $\epsilon\in\R_+$ and
suppose that $f$ has the property 
\begin{equation}\label{E: ExistsRight}
\st[f(r + dx)] = \st[f(r + dy)]\quad \text{in}\quad \R\,\cup\,\{\pm\infty\},
\end{equation} 
for all positive infinitesimals $d$x and $d$y. Then $\st(f(r + dx)) \in\R\,\cup\{\pm\infty\}$ 
is called the limit  of  $f$  as   $x$  approaches\,  $r$\, from the right, in symbols,
\begin{equation}\label{E: RightLimit}
\lim_{x\to r^+} f(x) = \st[f(r + dx)],
\end{equation}
where in the above formula $dx$ is an arbitrarily chosen positive infinitesimal. 
 
	{\bf (ii)} Let $f : X\to\R,\; X \subseteq\R$, be a real function and $(r -\epsilon, r) \subseteq X$ for some $r\in\R$ and some $\epsilon\in\R_+$.
Suppose, in addition, that    
\begin{equation}\label{E: ExistsLeft}
\st[f(r + dx)] = \st[f(r + dy)]\quad \text{in}\quad \R \cup\{\pm\infty\},
\end{equation} 
for all negative infinitesimals $dx, dy$.  Then $\st(f(r + dx)) \in \R \cup\{\pm\infty\}$ is called the limit of $f$ as $x$ approaches $r$  from
the left, in symbols, 
\begin{equation}\label{E: LeftLimit}
\lim_{x\to r^-} f(x) = \st(f(r + dx)),
\end{equation}
where $dx$ in the last formula is an arbitrarily chosen negative infinitesimal. \par
We refer to (\ref{E: RightLimit}) and (\ref{E: LeftLimit}) as our ``working formulas'' for the right and left-sided limits, respectively.
\end{definition}

\begin{remark}\label{R: Left and Right Limits}{\em (Left and Right Limits):} {\em  The comparison between the above definitions implies that 
$\lim_{x\to r} f(x)$ exists\, $if{f}$\, each of $\lim_{x\to r^+} f(x)$  and $\lim_{x\to r^-} f(x)$ exists and $\lim_{x\to r^+} f(x)=\lim_{x\to r^-}
f(x)$. In this case we have
\begin{equation}
				 \lim_{x\to r} f(x)=\lim_{x\to r^+}f(x)=\lim_{x\to r^-} f(x). 
\end{equation}
}\end{remark}

\begin{definition} \label{D: Limit at Infinity}{\em (Limit at Infinity):}
 
	{\bf (i)} Let $f : X \to \R,\; X \subset\R$, be a real function and $(a, \infty)\subset X$ for some $a\in\R$. Suppose that   
\begin{equation}\label{E: ExistsInfty}
\st(f(1/dx)) = \st(f(1/dy))\quad \text{in} \quad  \R\,\cup\{\pm\infty\},
\end{equation} 
for all positive infinitesimals $dx$  and $dy$. Then $\st(f(1/dx)) \in \R\,\cup\{\pm\infty\}$ is the limit of $f$  as $x$  goes to infinity, in
symbols,
\begin{equation}\label{E: WorkingFormulaInfty}
 \lim_{x\to \infty} f(x) = \st(f(1/dx)),
\end{equation}
where $dx$  is an arbitrarily chosen positive infinitesimal.

	{\bf (ii)} Let $f : X\to \R, X\subseteq\R$,\, be a real function and $(-\infty, b) \subset X$ for some $b\in\R$. Suppose that  
\begin{equation} \label{E: ExistsInfty-}
\st(f(1/dx)) = \st(f(1/dy))\quad \text{in}\quad \R\,\cup\{\pm\infty\},
\end{equation}
for all negative infinitesimals $dx$ and $dy$. Then $\st(f(1/dx))\in\R\,\cup\{\pm\infty\}$ is the limit of  $f$  as $x$ goes to minus infinity, in
symbols, 
\begin{equation} 	\label{E: WorkingFormulaInfty-}
\lim_{x\to -\infty}f(x) = \st(f(1/dx)),
\end{equation}
where $dx$ is an arbitrarily chosen negative infinitesimal. \par 
We refer to (\ref{E: WorkingFormulaInfty}) and (\ref{E: WorkingFormulaInfty-}) as our
``working formulas'' for the limits at infinity, respectively.
\end{definition}

\begin{remark} {\em (Existence of Limit at Infinity):} {\em As in the case of the usual limit (Remark~\ref{R: Existence of Limit}),  the condition
(\ref{E: ExistsInfty}) or (\ref{E: ExistsInfty-}) guarantees the existence of the limit value $\st(f(1/dx)) = L$ and the independence of the result on
the choice of $dx$. As before we should start with the calculation of  $L$ by (\ref{E: WorkingFormulaInfty}) or (\ref{E: WorkingFormulaInfty-}),
respectively, and when the calculations are over, we should check whether the value
$L$ depends on the choice of $dx$.  If, not, then $L$ is the correct answer. If $L$ depends on the choice of $dx$, then the corresponding limit
$\lim_{x\to\pm\infty}f(x)$ does not exist.
}\end{remark}
\begin{remark}{\em (Proper or Improper):} {\em If $\st(f(r + dx))$ or $\st(f(1/dx))$ is a real number, we say that the corresponding limits are
proper.  Otherwise (when a limit is $\infty$ or $-\infty$), we say that the limit is improper.  
}\end{remark}
\begin{remark}{\em (Unification):} {\em  We shall often unite  the right and left limits in the working formulae: 
\begin{equation}\label{E: WorkingFormulasUnified}
\lim_{x\to r^{\pm}}f(x) = \st(f(r + dx)),\quad\quad \lim_{x\to\pm\infty}f(x) = \st(f(1/dx)), 
\end{equation}
where in both formulae $dx$ is an arbitrarily chosen infinitesimal, positive or negative depending on the sign in $r^{\pm}$  or in $\pm\infty$,
respectively.   
}\end{remark}

\begin{theorem}{\em (A. Robinson):} The
above definitions of different type of limits are equivalent to the corresponding $\epsilon,\delta$-definitions. 
\end{theorem} 

\Proof: We refer the reader to Section~\ref{S: An Introduction to Infinitesimal Calculus} at the end of this paper. $\blacktriangle$
\section{Exercises on Limits}

	The main advantage of the formulas $\st(f(r+dx))$ and $\st(f(1/dx))$) over the standard $\epsilon,\delta$-definitions of 
$\lim_{x\to r^\pm}f(x)$ and $\lim_{x\to \pm\infty}f(x)$, respectively, is that $\st(f(r+dx))$ and $\st (f(1/dx))$ prescribe an algorithm for calculating 
the limits in terms of  $f$  and   $r$  only (without involvement of a candidate  $L$  for the limit value) : 
	
{\bf (a)} Evaluate  $f$  at the point  $r + dx$  (or at the point $1/dx$, respectively), where $dx$  is an infinitesimal, positive or negative, depending
on the sign in  $r^\pm$  or in $\pm\infty$, respectively. 

	{\bf (b)} Calculate the standard part of $f(r + dx)$ (or the standard part of $f(1/dx)$, respectively), following the rules in Section 1. 

	Let us write once again {\bf our working formulae:}
\begin{equation}\label{E: Working}
	\lim_{x\to r} f(x) = \st[f(r + dx)],
\end{equation}
where  $dx$ is an arbitrary non-zero infinitesimal, and 
\begin{align}
&\lim_{x\to r^\pm}f(x) = \st[f(r + dx)],  \label{E: WorkingBoth}\\
&\lim_{x\to \pm\infty}f(x) = \st[f(1/dx)],  \label{E: WorkingInftyBoth}
\end{align}
where in the last two formulae dx is an arbitrary infinitesimal, positive or negative, depending on the sign in $r^\pm$  or in  $\pm\infty$,
respectively. 

	Here are several examples of applications of our working formulae. The reader will observe that our method requires less sophistication in factoring
and less dependence on the Squeeze Theorem.

\begin{example}{\em $\lim_{x\to r} x^n = \st[f(r + dx)] = \st[(r + dx)^n] = (\st(r + dx))^n = r^n$. }\end{example}

\begin{example}{\em $\lim_{x\to \pm\infty}(1/x^n) =  \st[f(1/dx)] = \st(dx^n) =  (\st(dx))^n = 0$. }\end{example}

\begin{example}{\em $\lim_{x\to 0^\pm}(| x |/x) =  \st[f(dx)] = \st(|dx|/dx)  = \pm 1$, where  $dx$  is an arbitrary infinitesimal, positive or negative,
respectively. Notice that  $\lim_{x\to 0}(| x |/x)$ does not exist since $\lim_{x\to 0^+}(| x |/x)\not= \lim_{x\to 0^-}(| x |/x)$.}\end{example}

\begin{example} {\em $\lim_{x\to \pm\infty}(|x|/x) =\st[f(1/dx)] = \st(dx/|dx|)  = \pm 1$ (where, again, $dx$ is an arbitrary
infinitesimal, positive or negative, respectively). }\end{example}

\begin{example}   
   \begin{align*}
      \lim_{x\to (-1)^\pm}\frac{4x + 1}{x + 1}  
      &=  \st[f(-1 + dx)] 
      = \st\left(\frac{4(-1 + dx) + 1}{-1 + dx + 1}\right) \\
      &= \st\left(\frac{-3 + 4 dx}{dx}\right) 
      = \st(1/dx)\times\st(-3 + 4 dx) \\
      &= (\pm\infty)\times (-3) = \mp\infty,
   \end{align*}
{\em where $dx$ is an infinitesimal,  positive or negative, respectively. }\end{example}

\begin{example}\label{Ex: Factoring 1}{\em 
   \begin{align*}
      \lim_{x\to 1}\frac{x^3 + 4x^2 + x - 6}{x - 1}
      &=\st[f(1 + dx)]\\
	   &=  \st\left( \frac{(1+dx)^3 + 4(1+dx)^2 + 1 + dx - 6}{1 + dx - 1}\right)\\ 
      &=\st\left(\frac{12dx + 7dx^2 + dx^3}{dx}\right) \\
      &=  \st (12 + 7dx + dx^2)
      =  12
   \end{align*},
where $dx$ is a non-zero infinitesimal. These calculations will appeal to those students who prefer to expand the expressions  $(1 + dx)^3$ and 
$(1 + dx)^2$ (and collect the like-terms), rather than to factor the cubic function $x^3 + 4x^2 - x - 6$. Notice that the above calculations not only produce
the correct limit value 12, but also present a rigorous proof that 12 is, in fact, the  limit of the function, since the final result does not depend on the
choice of the infinitesimal  $dx$. }\end{example}
	
		Here is another example of an improper one-sided limit:

\begin{example}\label{Ex: Factoring 2}{\em  
   \begin{align*}
      \lim_{x\to 2^\pm}\left(\frac{x^3 - 9}{x^3 + x^2 - 7x + 2}\right)
      &=\st[f(2 + dx)] \\
      &= \st\left(\frac{(2+dx)^3 - 9}{(2+dx)^3 + (2+dx)^2 - 7(2+dx) + 2}\right)\\
	   &= \st\left(\frac{-1 + 12dx + 6dx^2 + dx^3}{9dx + 7dx^2 + dx^3}\right) \\
      &=  \st\left(\frac{-1 + 12dx + 6dx^2 + dx^3}{dx(9 + 7dx + dx^2)}\right) \\
		&= \st \left(\frac{1}{dx}\right)\times\st\left(\frac{-1 + 12dx + 6dx^2 + dx^3}{9 + 7dx + dx^2}\right)\\
		&= (\pm\infty)\times(-1/9) 
      =\mp\infty
   \end{align*}
   , where $dx$ is a non-zero infinitesimal, positive or negative, respectively. As in the previous
example,  these calculations will appeal to those students who prefer to expand the expressions  $(2 + dx)^3$ and $(2 + dx)^2$, rather than to factor the
cubic function $x^3 + x^2 - 7x + 2$.  }\end{example}

\begin{example} 
   \begin{align*}
      \lim_{x\to\pm\infty}\left(\frac{x}{\sqrt[4]{x^4+1}}\right) 
      &= \st[f(1/dx)] \\
      &= \st\left(\frac{1/dx}{\sqrt[4]{1/dx^4 +1}}\right)\\
      &=\st\left(\frac{|dx|}{dx}\times\frac{1}{\sqrt[4]{1 + dx^4}}\right)\\
      &= \st\left(\frac{|dx|}{dx}\right)\times\st\left(\frac{1}{\sqrt[4]{1 + dx^4}}\right) \\
      &= (\pm 1)\times\left(\frac{1}{\sqrt[4]{\st(1 + dx^4)}}\right) \\
      &= (\pm 1)\times 1 
      = \pm 1
   \end{align*},
{\em where  $dx$  is an infinitesimal,  positive or negative, respectively (notice that we use the fact that $\sqrt[4]{dx^4} = |dx|$). 
}\end{example}

\begin{example}\label{Ex: Sine} {\em $\lim_{x\to\pm\infty}[\sin(x)/x] = \st [f(1/dx)] = \st(dx)\times\sin(1/dx)] = 
\st(dx)\times\st[\sin(1/dx)] = 0 \times\st[sin(1/dx) = 0$, since $\st[\sin(1/dx)]$ is a well defined (although explicitly unknown) real number in
[-1, 1]  (see Remark~\ref{R: TwoPrinciples} about the range of\, $\sin x$). 
}\end{example}

\begin{remark}\label{R: Squeeze Theorem} {\em (No Need of the Squeeze Theorem):} {\em We believe that the calculations in the last example, based on
the product rule for standard part \st\, (Theorem~\ref{T: Properties of st}), are shorter and simpler than the usual Squeeze Theorem arguments.
Recall that the product formula for limits is non-applicable in the case of the function\, $\sin{x}/x$\, since $\lim_{x \to\infty}\sin{x}$ does not
exist in \R\, (even as an improper limit). In contrast, we have $\st[dx \sin(1/dx)] = \st(dx) \st[\sin(1/dx)]$ since each of $\st(dx)$ and
$\st[\sin(1/dx)]$ exists in \R, by Theorem~\ref{T: Asymptotoc Expansion} and Lemma~\ref{L: Uniqueness}. In general, the infinitesimal method for
calculating limits is less dependent on the Squeeze Theorem, and more rarely requires the use of inequalities, compared with the usual Weierstrass
$\epsilon,\delta$-method. We consider that feature of the Infinitesimal Calculus as an advantage over the conventional standard methods, taking
into account that the students in calculus are rarely in love with inequalities. }\end{remark} 

	Although our text is devoted to limits only, we shall shortly mention the definitions of derivative  and differential  in terms of infinitesimals:

\begin{remark}{\em (Derivative):} {\em  We define the derivative  by 
   \[f'(x) = \st\left(\frac{f(x + dx) - f(x)}{dx}\right),\]
where $dx$ is an arbitrary non-zero infinitesimal. For example, 
\begin{align*}
   \left(x^3\right)' 
   &= \st\left(\frac{\left(x + dx\right)^3 - x^3}{dx}\right) 
   = \st\left(\frac{x^3 + 3x^2 dx + 3 x dx^2 + dx^3 - x^3}{dx}\right) \\
   &=  \st\left(\frac{3 x^2 dx + 3 x dx^2 + dx^3}{dx}\right) 
   = \st\left(3 x^2 + 3 x dx + dx^2\right) = 3 x^2.
\end{align*}
Similarly, we define differential by $dy =
f'(x) dx$, where $dx$ is an infinitesimal. Thus, the Leibniz notation  $dy/dx = f'(x)$ holds ``by the definition'' of $dy$  for all non-zero
infinitesimals $dx$.   } \end{remark}

\noindent {\bf Summary:} We summarize the properties of our working formulae~(\ref{E: Working})-(\ref{E: WorkingInftyBoth}): 

	{\bf 1.}  Our working formulae, if considered as definitions, are equivalent of the usual $\epsilon,\delta$-definitions of the corresponding limits at any
level of generality and rigor (see the next section). 

	{\bf 2.} Our working formulae are free of a candidate, $L$, for the limit value. Hence we do not need to guess and check (because there is
nothing to guess). Also, our working formulae ``really work'' in the sense that they produces the correct value of $L$, as demonstrated by the
above examples. 

	{\bf 3.} Our method requires less sophistication in factoring (Examples~(\ref{Ex: Factoring 1})-(\ref{Ex: Factoring 2}) and it is less dependent
on the Squeeze Theorem compared with the conventional methods (Example~(\ref{Ex: Sine}) and Remark~(\ref{R: Squeeze Theorem})).

	{\bf 4.} Under the assumption that ``the limit exists'', our working formulae~(\ref{E: Working})-(\ref{E: WorkingInftyBoth}) are free of
quantifiers, as opposed to the three non-commuting quantifiers ``$\forall, \exists, \forall$'' in the $\epsilon, \delta$-definition of limit. On the other
hand, each of the  criterions for existence of limit ((\ref{E: Exists}), (\ref{E: ExistsRight}), (\ref{E: ExistsLeft}), (\ref{E: ExistsInfty}), (\ref{E:
ExistsInfty-})) involves two commuting quantifiers ``$\forall, \forall$'' only,  as opposed to the four non-commuting quantifiers ``$\forall, \exists,
\forall,\forall$'' in its standard counterpart. Thus, our method is easier to apply to rigorous proofs when (and if) the teacher decides to do rigorous
proofs. In fact, in our method, the rigorous proof that ``$L$  is, actually, the correct limit value'' coincides with the ``calculation of  $L$.''
	
	As a result, the Calculus becomes  to ``calculate''  again, as it was originally designed to do. 

\begin{remark} {\em (Infinitesimals in Mathematica):} {\em Assume that you are already ``addicted'' to calculators and computers and that you are
not planning to ``quit'' any time soon. Assume that your students have already purchased expensive calculators and are eager to calculate ``anything
which comes along.'' Assume that your university has already spent a lot of money on buying computers and the spending has to be somehow justified.
And to complete the scenario assume, finally, that you have just won a generous grant from NSF for ``using technology in teaching calculus.''
Under these circumstances your dilemma will be ``how to reconcile the computers with infinitesimals''? The good news is that computers are able to
handle infinitesimals, and they actually work with infinitesimals anyway. Take, for example, ``Mathematica''. Have you ever thought about how
Mathematica calculates limits ? It might occur to you that the computer evaluates the function  $f$ at finitely many points and announces one these
values for the ``correct answer'' ? Or, perhaps, the computer has simply memorized the limits of all possible functions ``you will ever ask it
for''? The answer is ``neither of the above.'' Rather, Mathematica calculates the formal (Taylor) asymptotic expansion of  $f(r + dx)$ by the
command ``Series'', treating  $dx$  as a ``formal variable'' and truncates all terms in the series but the first by the command ``$/{\bf .}\,
dx\to0$". This procedure has very little to do with  the $\epsilon,\delta$-definition of limit and it is almost identical to the operation
``taking the standard part'' discussed in this article. The framework of these calculations is the field $\R((dx))$ of formal Laurent series with
real coefficients and formal variable, denoted  by $dx$. Notice that the field $\R((dx))$ is non-Archimedian and the formal
variable, $dx$, if considered as an element of $\R((dx))$, is a positive infinitesimal. So what ? Well, it means that Mathematica (believe it
or not) calculates limits through infinitesimals in the framework of the field $\R((dx))$. The formula $\st[f(x + dx)] = Series[f(x + dx),{dx,0,1}]
/. dx\to 0$  can be used for calculating the ``standard part'' of $f(x + dx)$ in Mathematica if you decide to do so. The author of this article,
however, is unable to see any pedagogical merits of this activity unless, perhaps, for the purpose of a better understanding of how Mathematica
works.  }\end{remark}


\section{An Introduction to Infinitesimal Calculus}\label{S: An Introduction to Infinitesimal Calculus}

	Here we present a short introduction to the modern Infinitesimal Calculus known as well as A. Robinson's Nonstandard Analysis. We would like to 
assure the reader that the usual background in real analysis is more than enough to follow this text. For more detailed exposition we shall 
refer to (H. J. Keisler~\cite{hjKeisE}-\cite{hjKeisF}) and (T. Lindstr{$\emptyset$}m~\cite{tLind}), where the reader will find more references to the
subject. For a really short  (although somewhat dense) exposition of both axiomatic and sequential approaches to nonstandard analysis we refer to T.
Todorov~\cite{tTod},  p. 685-688. We shall restrict our exposition to the nonstandard treatment of proper limits of the form $\lim_{x\to r}f(x)$ only,
and leave the improper limits, as well as the limits at  infinity, to the reader. All results in Section 2, follow as particular  cases. 

		Although the nonstandard analysis arose historically in a close connection with model theory and mathematical logic (A. Robinson~\cite{aRob66}), it
is completely possible to construct it in the framework of the standard analysis, i.e. assuming only the properties of the real numbers (along with
the Axiom of Choice). The method (due to W. A. J. Luxemburg), is known as the ultrapower construction or constructive nonstandard analysis:
	
	{\bf 1.} Let \N\, be the set of the natural numbers and $\mathcal{P}(\N)$ be the power set of \N. Let  $\mu: P(\N)\to \{0, 1\}$ be a two-valued finitely
additive measure such that $\mu(A) = 0$ for all finite $A\subset\N$\, and $\mu(\N) = 1$. We shall keep  $\mu$  fixed in what follows.	 

\begin{remark}{\em (Existence of $\mu$):} {\em To show that there exists a measure with these properties, it suffices to take a free ultrafilter $\mathcal{U}$
on
\N\, and define $\mu$ by $\mu(A) = 1$ for $A\in\mathcal{U}$ and $\mu(A) = 0$ for $A \notin \mathcal{U}$. Recall that a non-empty set $\mathcal{U}$ of
subsets of
\N\, is called a free ultrafilter on \N\, if it satisfies the following four properties: {\bf (a)} $\mathcal{U}$ is closed under
intersection;\; {\bf (b)} If $A, B \subseteq \N$, then $\mathcal{U}\ni A\subseteq B$ implies $B\in \mathcal{U}$;\; {\bf (c)} For any $A\subseteq\N$
exactly one of the following is true: $A\in \mathcal{U}$ or $\N\setminus A\in \mathcal{U}$;\;	{\bf (d)} $\bigcap_{A\in\mathcal{U}}A
= \emptyset$. Recall that the existence of free ultrafilters on
\N, follows from the Axiom of Choice (H. J. Keisler \cite{hjKeisF}, p. 49). We should mention
that the familiarity with the theory of ultrafilters is not necessary for the understanding of what follows. }\end{remark}


	The next properties of\, $\mu$\,  follow immediately from the definition:

\begin{lemma}{\em (Properties of $\mu$):}\label{L: Properties of mu}  Let $A, B \subset\N$. Then: 

	{\bf (a)} $\mu(A\cup B) = 1$ $\iff$ $[\,\mu(A) = 1$ or $\mu(B) = 1\,]$. In particular, for any $A\subseteq\N$ exactly one of $\mu(A) = 1$ and
$\mu(\N\,\setminus A) = 1$ is true.

	{\bf (b)} $\mu(A) = 1$ for all co-finite sets $A$ of \N. In particular, $\mu(\N)=1$.


	{\bf (c)} $\mu(A\cup B) = 0 \iff [\mu(A) = 0$ and $\mu(B) = 0]$. 
 
	{\bf (d)} $\mu(A)=\mu(B) = 1 \iff \mu(A\cap B) = 1$.

	{\bf (e)} $A\subseteq B \subseteq\N$ and $\mu(A) = 1$ implies $\mu(B) = 1$.
\end{lemma}

{\bf 2.} Let $\R^{\N}$ be the set of all sequences of real numbers considered as a ring under the usual pointwise operations. Define an equivalence
relation $\sim$ in $\R^{\N}$ by:  $(a_n) \sim (b_n)$ if $a_n = b_n$  $a. e.$ (where ``$a. e.$'' stands for ``almost everywhere''), 
i.e. if $\mu(\{\,n\, \mid\, a_n = b_n\,\}) = 1$. Then the factor space 
$\starR =\R^{\N}/\sim$ defines a set of nonstandard real numbers (or hyperreals). We shall denote by $\langle a_n\rangle$ the equivalence class
determined by the sequence $(a_n)$. We also define the embedding
$\R\,\subset\starR$ by $r\to \langle r, r, r,\dots\rangle$. In what follows we shall identify notationally a given real numbers $r$ with its image in
\starR. The addition and multiplication in
\starR\, is inherited from $\R^{\N}$. The order relation in \starR\, is defined by: $\langle a_n\rangle\leq \langle b_n\rangle$ if  $a_n\leq b_n$ 
holds 
$a. e.$,  i.e. if 
$\mu(\{\, n\,
\mid\, a_n
\leq b_n\,\}) = 1$. We define also $|x| = \max\,\{x, -x\}$. Notice that we define one specific nonstandard extension  \starR\, of the reals \R\, which
depends, in general, on the choice of the measure $\mu$. We should mention that the different fields of the form \starR\, (corresponding to different
measures $\mu$) are not necessarily isomorphic to each other. We also have $\text{card}(\!\starR) = \text{card}(\!\R)$.

\begin{theorem} {\em \starR}\, is a totally ordered non-Archimedean field  containing \R\, as a totally ordered subfield. \end{theorem}
\Proof \starR\, is a ring since $\R^{\N}$ is a ring. To show that \starR\, has no zero divisors, assume that $\langle a_n\rangle\langle 
b_n\rangle = 0$ in \starR, i.e. $\mu(\{\, n\, \mid\ a_n\, b_n = 0\,\}) = 1$. Denote $A = \{\,n\, \mid\, a_n = 0\,\}$ and 
$B = \{\, n\, \mid b_n = 0\,\}$ and observe that $\{\,n \mid a_n\,b_n = 0\,\} = A \cup B$, since \R\, has no zero divisors. Hence either $\mu(A) = 1$  or
$\mu(B) = 1$, by Lemma~\ref{L: Properties of mu}, i.e. either $\langle a_n\rangle = 0$, or $\langle b_n\rangle = 0$, as required. 
To show that the non-zero
elements in \starR\, are multiplicative invertible, assume that 
$\langle a_n\rangle \not= 0$ in  \starR, i.e. $\mu(\{\, n\, \mid\, a_n \not=  0\}) = 1$. 
Denote $\{\, n\mid a_n \not= 0\,\} = C$ and define $(b_n)\in \R^{\N}$ by
$b_n = 1/a_n$  if $n \in C$ and anyhow (say,  $b_n = 1$)  if $n\in\N\setminus C$. We have 
$C\subseteq \{\, n\, \mid  a_n\, b_n  = 1\,\}$ which implies $\mu(\{\, n\, \mid a_n b_n  = 1\,\}) = 1$, by Lemma~\ref{L:
Properties of mu}. Thus $\langle a_n \rangle\langle b_n\rangle = 1$, as required. To show the trichotomy of the order relation, assume that $\langle
a_n\rangle\not=\langle b_n\rangle$ and denote $A = \{\, n\, \mid  a_n\leq b_n\,
\}$. We have $\N\setminus A = \{\, n\, \mid a_n > b_n\,\}$ and, thus, exactly one of $\mu(A) = 1$ or $\mu(\N\setminus A) = 1$ is true, 
by Lemma~\ref{L:
Properties of mu}. That is $\langle a_n\rangle \leq \langle b_n\rangle$ or  $\langle a_n\rangle > \langle b_n\rangle$ which is equivalent to 
$\langle
a_n\rangle < \langle b_n\rangle$ or  $\langle a_n\rangle > \langle b_n\rangle$, as required, since $\langle a_n\rangle \not= \langle b_n\rangle$, by
assumption. The rest of the properties of the totally ordered field can be proved similarly. The embedding $\R\,\subset \starR$  is obviously field and
order preserving. To show that  \starR\, is
non-Archimedean, observe that
$ m <
\langle n\rangle$  in
\starR\, for any m  in
\N,  (where  m  is considered as an element of \starR) since the set $\{\, n\, \mid m  <  n\, \}$ is co-finite, and hence, of measure $1$. $\blacktriangle$


\begin{example} {\em $\langle 1/n\rangle, \langle1/n^2\rangle, \langle 1/\ln{n}\rangle, \langle e^{-n}\rangle$ are positive infinitesimals (different from each
other) and $\langle n\rangle, \langle n^2\rangle, \langle \ln{n}\rangle, \langle e^n\rangle$ are positive infinitely large numbers (also different
from each other). The number $\langle 3+ 1/n\rangle = 3 +\langle 1/n\rangle$ is finite
(but not real). Let us take the first example: for any $m\in \N$,  the set $\{\, n\, \mid  0 < 1/n < 1/m\,\}$ is co-finite, hence, of measure 1. Therefore,
$0 < \langle 1/n\rangle < 1/m$  in \starR, i.e. $\langle 1/n\rangle$ is a positive infinitesimal. The rest of the examples are treated similarly. 
}\end{example}

	{\bf 3.} It is clear that $\R\subset \mathcal{F}(\!\starR),  \mathcal{I}(\!\starR) \subset  \mathcal{F}(\!\starR),  \R\cap  \mathcal{I}(\!\starR) =
\{0\},
\mathcal{F}(\!\starR) \cap
\mathcal{L}(\!\starR)  = 
\emptyset$  and 
$\mathcal{F}(\!\starR) \cup\mathcal{L}(\!\starR) =  \starR$.  From the above definition it follows easily that $\mathcal{F}(\!\starR)$ is a totally
ordered integral domain and
$\mathcal{I}(\!\starR)$  is a convex maximal ideal in
$\mathcal{F}(\!\starR)$. Hence, $\mathcal{F}(\!\starR)/\mathcal{I}(\!\starR)$ is a totally ordered field which is isomorphic to \R\, as totally ordered
fields. The canonical homomorphism
$\st:\mathcal{F}(\!\starR)
\to \R$ is called the {\em standard part mapping}. Notice that  $\st\langle a_n\rangle$  exists for any bounded sequence $(a_n)$ in $\R^{\N}$ and $\st\langle
a_n\rangle =
\lim_{n\to\infty} a_n$  for any convergent $(a_n)$. Conversely, if $\langle a_n\rangle$ is a finite number, then $\st\langle a_n\rangle = \lim_{n\to\infty}
a_{k_n}$  for some subsequence
$(a_{k_n})$  of  $(a_n)$  such that $\mu(\{\, k_n\, \mid n\in\N\,\}) = 1$. The following result follows immediately: 

\begin{theorem}\label{T: Standard Part Mapping}

	{\bf (i)} Let  {\em $x\in\starR$}. Then {\em $x\in\mathcal{F}(\!\starR)\quad if{f}\quad  x = r + dx$}  for some $x\in\R$ and some
{\em $dx\in\mathcal{I}(\!\!\starR)$}. 
	
	{\bf (ii)} If {\em $x\in\mathcal{F}(\!\starR)$}, then the presentation  $x  =  r + dx$ is unique and  $r = \st(x)$. In particular,  $\st(r) = r$  for any
$r\in\R$. 

	{\bf (iii)} The standard part mapping is order preserving in the sense that  $x\leq y$ in {\em $\mathcal{F}(\!\!\starR)$}  implies  $\st( x) \leq \st
(y)$ in
\R. 
\end{theorem}

	{\bf 4.} Let $X\subseteq\R$. Then the set\,$\asterisk X = \{\,\langle x_n\rangle \in\starR\, \mid x_n\in X\;  \text{a. e.}\,\}$ is called the
nonstandard extension of  $X$.  For any Ê$X \subseteq\R$ we have $X\subseteq\asterisk X$ and $X = \asterisk X\quad if{f}\quad X$ is a
finite set. The above definition holds also in the case when $X\subseteq\R^d\; (d\in\N$). If $X\subseteq\R$ and $Y\subseteq\R$, then $\asterisk(X\times Y)
= \asterisk X\times\asterisk Y$. In particular, we have $\asterisk(\!{\R}^d) = (\!\starR)^d$, so we can write simply  $\asterisk\R^d$.


\begin{example} {\em Let $\Q, \Z,\N, \mathbb{E},\mathbb{O}, \mathbb{P}$, etc. be the sets of the rational, integer, natural, even, odd, prime, etc. numbers,
respectively. Then the elements of $\asterisk\Q, \asterisk\Z, \asterisk\N, \asterisk\,\mathbb{E}, \asterisk\,\mathbb{O}, \asterisk\,\mathbb{P}$, etc. will
be called nonstandard rational numbers, nonstandard integer numbers (hyperintegers), nonstandard natural numbers (hypernatural numbers), etc., respectively. The
set of the {\em infinitely large natural numbers}, i.e. the infinitely large numbers in \asterisk\N, will be denoted by $\N_\infty$. We have  $\asterisk\N =
\N\,\cup\N_\infty$.  If $a, b \in \R$,  then 
\[
\asterisk\,\{\,x\in\R\,\mid  a\leq x \leq b\,\} = \{\,x \in \starR\, \mid  a\leq x \leq b\,\},
\]
which will be denoted for short by $\asterisk\,[a, b]$. It follows that $r + dx\in \asterisk\,[a, b]$ for
all $r \in [a, b)$ and all non-negative infinitesimal $dx$. Similarly, we have  
\[
\asterisk\,\{\,x \in \R\, \mid  a < x < b\,\} = \{\, x\in\starR\, \mid  a < x < b\,\},
\]
which will be denoted for short by $\asterisk\,(a, b)$. We have $r + dx \in \asterisk\,(a, b)$ for all $r \in (a, b)$ and all infinitesimal  $dx$ .
}\end{example}

\begin{theorem}\label{T: Adherent Point} {\em (Adherent Point):} Let $r\in\R$ and $X\subseteq\R$. Then $r$ is a non-trivial adherent point  of $X\; if{f}$\, there
exists {\em
$dx\in\starR$} such that $dx\not= 0, dx\approx0$ and {\em $r + dx\in\asterisk X$} (or, equivalently,\;  $if{f}$\,  there exists {\em $x \in \asterisk X$} such
that  $r\not= x$ and $r \approx x$). 
\end{theorem}

\Proof: $(\Rightarrow)$ For any $n\in\N$ the set X$_n = \{\, x\in X \mid 0 <  | x - r | <  \frac{1}{n}\,\}$ is non-empty, by assumption. Hence, by the Axiom of
Choice, there exists $(x_n)$ in $\R^{\N}$ such that $x_n\in X_n$  for all $n\in\N$. Now, $dx = \langle x_n\rangle - r$ is the infinitesimal we are looking
for. Indeed, we have $0  <  |dx|  < \langle\frac{1}{n}\rangle$, hence $dx\not= 0$ and $dx\approx0$,  since $\langle\frac{1}{n}\rangle\approx 0$.  Also  $r +
dx =
\langle x_n\rangle\in \asterisk X$, since $x_n\in X_n$  for all $n \in\N$.
	$(\Leftarrow)$ We have $dx\not= 0, dx\approx 0, r + dx\in\asterisk X$ for some $dx\in\starR$, by assumption. Suppose $m\in\N$ and observe that $0 < |dx| <
1/m$. We have $dx = \langle \epsilon_n\rangle$ for some $(\epsilon_n)$ in
$\R^{\N}$. The set  $\{\, n \mid   0 <  |\epsilon_n| < 1/m,\;  r + \epsilon_n \in X\,\}$ is of measure $1$, hence, it is non-empty. The
latter means that $r$  is a non-trivial adherent point of  $X$. $\blacktriangle$ 
	
{\bf 5.} Let $f : X \to \R$ be a real function, where $X\subseteq\R$. Then the function  $\asterisk f : \asterisk X \to \starR$,  defined by $\asterisk
f(\langle x_n\rangle) = \langle f(x_n)\rangle$ for all $\langle x_n\rangle \in \asterisk X$, is called the nonstandard extension of $f$  since $\asterisk
f(r) = f(r)$ for all $r\in X$. The above definition holds also in the case  $X\subseteq \R^d\; (d\in\N$). 

\begin{theorem}\label{T: A. Robinson} {\em (A. Robinson):} Let $r$ be a non-trivial adherent point of \par 
\noindent $f: X\to\R,\, X\subseteq\R$, be a real function
and
$L\in\R$. Then
$\lim_{x\to r} f(x)  =  L\quad if{f}$\quad {\em $\asterisk f(r + dx) \approx L$} for all {\em $dx\in\starR$} such that $dx\not= 0,\; dx\approx 0$ and {\em $r
+ dx
\in
\asterisk X$}. If the limit exists in \R, then {\em $\lim_{x\to r} f(x) = \st(\asterisk f(r + dx))$}. 
\end{theorem}

\Proof: ($\Rightarrow$)  Let $\epsilon\in\R_+$. By assumption, there exists $\delta\in
R_+$ such that for all $x\in X,\quad 0 < | x - r | < \delta$ implies $| f(x) - L | < \epsilon$. Let $dx\not= 0,\; dx\approx 0$ and $r + dx \in
\asterisk X$ for some $dx\in\starR$. Notice that  $dx$  exists, by Theorem~\ref{T: Adherent Point}, since\, $r$\, is a non-trivial adherent point of
$X$, by assumption. We have  $dx = \langle \epsilon_n\rangle$ for some sequence $(\epsilon_n)$ in $\R^{\N}$. Next, we define the sets: 
\[
A_\delta = \{\, n\, \mid  0 <  | \epsilon_n |  <  \delta\; \text{and}\; r + \epsilon_n  \in X\, \}\; \text{and}\; B_\epsilon = \{\, n\, \mid\, | f(r +
\epsilon_n) - L |  <
\epsilon\,\}.
\]
We have $\mu(B_\epsilon) = 1$, by Lemma~\ref{L: Properties of mu}, since $A_\delta\subseteq B_\epsilon$ and $\mu(A_\delta) = 1$, by assumption.
Recapitulating, we have
$|\asterisk f(r + dx) - L | < \epsilon$  for all  $\epsilon \in
\R_+$, which means that $\asterisk f(r + dx) \approx  L$,  as required. 

	($\Leftarrow$) Assume (on the contrary) that  $\lim _{x\to r} f(x) = L$ is false. Thus, there exists $\epsilon\in \R_+$ such that the sets 
\[
	X_n  =  \{\, x \in X\, \mid\,   0 < | x - r | <  \frac{1}{n}\; \text{and}\; | f(x) - L | > \epsilon\,\},
\]
are non-empty for each $n \in\N$. By the Axiom of Choice, there exists a sequence $(x_n)$ in $\R^{\N}$ such that $x_n\in X_n$ for all $n\in \N$. Define
$dx\in\starR$  by $dx = \langle x_n\rangle - r$.  We have $0  <  | dx |  <  \langle\frac{1}{n}\rangle$ hence $dx \not= 0$ and $dx \approx 0$.  Also  $r
+ dx =  
\langle x_n\rangle\in \asterisk X$  and $| \asterisk f(r + dx) - L | > \epsilon$.  The latter means that $\asterisk f(r + dx)- L$ is a non-infinitesimal, a
contradiction.   The formula $\lim_{x\to r} f(x) = \st(\asterisk f(r + dx))$ follows directly from $\asterisk f(r + dx) \approx L$ after applying the
standard part mapping to both sides and taking into account that  $\st(L) = L$ since $L$ is a real number. $\blacktriangle$


	In order to eliminate completely the ``candidate'' for limit $L$ from our theory, we have to present a nonstandard characterization of the existence of a proper
limit in terms of infinitesimals similar to the Cauchy convergence criterion: 

\begin{theorem}\label{T: Existenced}{\em (Existence):} Let  $f :  X \to \R,\,  X\subseteq \R$,  be a real function and let  $r\in
\R$  be a non-trivial adherent point of  $X$. Then the following are equivalent: 
	
	{\bf (i)}	The limit  $\lim_{x\to r} f(x)$ exists in \R.

	{\bf (ii)} 	$f$ is fundamental  (or Cauchy) toward  $ r$  in the sense that 
\[
		(\forall\epsilon\in \R_+)(\exists\delta\in \R_+)(\forall x, y \in X ) [ 0 <  |x - r|,\,  |y - r|  <  \delta \Rightarrow | f(x) - f(y) |  < \epsilon ].
\]

	{\bf (iii)} {\em $\asterisk f(x) \approx \asterisk g(y)$} for all {\em $x, y\in \asterisk X$} such that  $x\not= r, y \not= r,  x \approx r$ and $y
\approx r$. 

	{\bf (iv)} {\em $\asterisk f(r+dx) \approx \asterisk g(r+dy)$}  for all  $dx, dy\not= 0, dx, dy \approx 0$ such that {\em $r + dx,\; r + dy\in\asterisk
X$}. 

	{\bf (v)} {\em $\st(\asterisk f(r + dx)) = \st(\asterisk f(r + dy)) \in \R$} (but never become $\pm\infty$) for all  $dx, dy\not= 0, dx, dy \approx 0$ such
that {\em $r + dx,\; r + dy \in \asterisk X$}. 

	{\bf (vi)} {\em $(\exists h\in\mathcal{I}(\!\starR_+))(\forall dx, dy \in\starR)$}
{\em \[
\left [r + dx, r + dy\in \asterisk X\; \text{and}\; 0 < |dx|, |dy| <  h \right] \implies \left[\asterisk f(r + dx) \approx \asterisk f(r +
dy)\right],
\]}
where {\em $\mathcal{I}(\!\starR_+)$} denotes the set of the positive infinitesimals in {\em \starR}.
\end{theorem}
\Proof: (i) $\Leftrightarrow$ (ii) is the Cauchy Criterion for existence of limits (Alan F. Beardon~\cite{aBear}, Theorem 4.4.1, p. 57). 
		
	(i) $\Rightarrow$(iii): We have $\asterisk f(x) \approx L$  and  $\asterisk f(y) \approx L$  for the same $L\in\R$, by Theorem~\ref{T: A. Robinson}, hence,
$\asterisk f(x) \approx \asterisk f(y)$, as required. 

	(iii)$\Leftrightarrow$(iv) follows immediately by letting $x = r + dx$ and $y = r + dy$. 

	(iv)$\Rightarrow$ (vi) in a trivial way.

	(vi)$\Rightarrow$(i) : We have  $h = \langle h_n\rangle$ for some sequence $(h_n)$ in $\R^{\N}$. Without loss of generality we can assume that  $h_n > 0$
for all $n
\in\N$. Now, suppose (for contradiction) that (i) fails, i.e. there exists  $\epsilon\in\R_+$  such that  $A_n\not=\emptyset$  for all $n\in\N$, where 
\[
A_n = \{\,\langle x, y\rangle\in X\times X\,  \mid\;  0 < | x - r |,\; | y - r | <  h_n\;  \text{and}\;  | f(x) - f(y) | \geq \epsilon\,\}.
\]
 Hence (by Axiom of Choice), there exists a sequence $(x_n, y_n)$ in 
$X^{\N} \times X^{\N}$ such that  $(x_n, y_n) \in A_n$ for all $n\in\N$. We define the nonstandard numbers $\langle x_n\rangle, \langle
y_n\rangle \in
\asterisk X$ and observe that 
\[
0 < | \langle x_n\rangle - r |,\, |\langle y_n\rangle - r | < \langle h_n\rangle\quad \text{and}\quad |\asterisk f\langle x_n\rangle) - \asterisk f(\langle
y_n\rangle) |
\geq  \epsilon,
\]
in \starR\,, by the choice of $(x_n)$ and $(y_n)$. Thus, $\asterisk f( \langle x_n\rangle)-\asterisk f(\langle y_n\rangle)$ is a non-infinitesimal, which
contradicts (vi) for $ h = \langle h_n\rangle,\,  dx =\langle x_n\rangle - r$ and  $dy =  \langle y_n\rangle- r$. 
	
(i)$\Rightarrow$(v) : We have $\st(\asterisk f(r + dx)) = L$  and  $\st(\asterisk f(r + dy)) = L$  for the same $L\in\R$, by Theorem~\ref{T: A. Robinson}, hence, 
$\st(\asterisk f(x)) =  \st(\asterisk f(y)) \in\R$, as required. 
	
(iv) $\Rightarrow$ (v) : $\st(\asterisk f(r + dx)) = \st(\asterisk f(r + dy)) \in\R$ implies, in particular, that $\asterisk f(r + dx)$ and $\asterisk f(r + dy)$
are finite numbers, thus, $\asterisk f(x) \approx \asterisk f(y)$ follows. $\blacktriangle$

\noindent {\bf Simplified Notation:} For the purpose of teaching and explicit calculations we recommend the following simplified notations:
	
	{\bf (a)} We shall skip the asterisks in front of $\asterisk f$, writing  simply $f$. This is perfectly justifiable since $\asterisk f$ is an extension of
$f$. 
	
	{\bf (b)} If  $X\subseteq \R$, then we shall sometimes write simply $X$ meaning $\asterisk X$. For example, we shall write $[a, b]$ meaning, actually, 
\[
\asterisk [a, b] = \{\,x \in\starR\, \mid a \leq x \leq b\,\}.
\]
	
	{\bf (c)} Finally, we shall write $(-\infty,\,\infty)$ for both \R\,  and \starR\, leaving the reader to figure out from the context which one we mean. 

	{\bf (d)} We prefer to use the terminology ``hyperreal numbers'' rather than ``nonstandard numbers'' (to avoid the shocking effect of the word
``nonstandard").
	
	{\bf (e)} We preserve our rights to come back to the more precise $\asterisk$-notation when (and if) needed. \\


	{\bf Acknowledgment:} The author thanks Ivan Penkov, Richard Pollard, Robert Wolf and the referee of this paper for
the many corrections and useful suggestions to the earlier versions of the manuscript.
\newpage


\end{document}